\documentclass[article]{siamart}

\usepackage{amsmath,mathtools}
\usepackage{latexsym}
\usepackage{graphics}
\usepackage{algorithm}
\usepackage{algorithmicx}
\usepackage{graphicx}
\usepackage{booktabs}
\usepackage{epstopdf}
\usepackage{microtype}
\usepackage[english]{babel}
\usepackage{amssymb}
\usepackage{color}
\usepackage{verbatim}
\usepackage{geometry}
\usepackage{marginnote}
\usepackage{autonum}
\usepackage{caption}
\usepackage{wrapfig}
\usepackage{sidecap}
\usepackage[font=scriptsize]{caption}
\usepackage{subcaption}
\usepackage{hyperref}
\usepackage{makecell}
\usepackage{algpseudocode}
\usepackage{xcolor}

\algnewcommand{\LineComment}[1]{\State \(\triangleright\) #1}
\graphicspath{ {YAF-Figures/FinalFigsMSD/} }

\renewcommand{\theequation}{\arabic{section}.\arabic{equation}}

\newcommand{\beq}{\begin{eqnarray}}
\newcommand{\eeq}{\end{eqnarray}}
\newcommand{\be}{\begin{equation}}
\algnewcommand\algorithmicinput{\textbf{Input:}}
\algnewcommand\INPUT{\item[\algorithmicinput]}
\algnewcommand\algorithmicoutput{\textbf{Output:}}
\algnewcommand\OUTPUT{\item[\algorithmicoutput]}
\newcommand{\ee}{\end{equation}}
\newcommand{\beqq}{\begin{eqnarray*}}
\newcommand{\eeqq}{\end{eqnarray*}}

\crefname{enumi}{assumption}{assumptions}
\headers{Regularization with the Picard parameter}{Eitan Levin and Alexander Y. Meltzer}
\title{Estimation of the Regularization Parameter in Linear Discrete Ill-Posed Problems Using the Picard parameter\thanks{This research was supported by the Israel Science Foundation (grant No. 132/14) and by Rosa and Emilio Segr\'{e} Research Award. }}

\author{
  Eitan Levin\thanks{Department of Condensed Matter Physics, Weizmann Institute of Science, 76100 Rehovot, Israel \newline
    (\email{eitan.levin@weizmann.ac.il}, \email{alexander.meltzer@weizmann.ac.il}.)}
  \and
  Alexander Y. Meltzer\footnotemark[2]
}
\begin{document}
\maketitle
\begin{abstract}
 Accurate determination of the regularization parameter in inverse problems still represents an analytical challenge, owing mainly to the considerable difficulty to separate the unknown noise from the signal. We present a new approach for determining the parameter for the general-form Tikhonov regularization of linear ill-posed problems. In our approach the parameter is found by approximate minimization of the distance between the unknown noiseless data and the data reconstructed from the regularized solution. We approximate this distance by employing the Picard parameter to separate the noise from the data in the coordinate system of the generalized SVD. A simple and reliable algorithm for the estimation of the Picard parameter enables accurate implementation of the above procedure. We demonstrate the effectiveness of our method on several numerical examples\footnote{A MATLAB-based implementation of the proposed algorithms can be found at \href{https://www.weizmann.ac.il/condmat/superc/software/}{https://www.weizmann.ac.il/condmat/superc/software/}}.
\end{abstract}



\section{Introduction}\label{sec:intro}

The Tikhonov regularization method \cite{Tikhonov1978} is one of the most widely applied methods for solution of linear ill-posed problems. It is well known that the accuracy of the solution obtained using the Tikhonov regularization method depends crucially on the chosen regularization parameter. This parameter is often obtained using either one of the following methods - the Generalized Cross-Validation (GCV) \cite{splines,GCV2}, L-curve \cite{LCurveAnal,LCurve}, Quasi-optimality \cite{QO,QO1}, Stein's Unbiased Risk Estimate (SURE) \cite{SURE,blackBoxSure}, or other methods. However, none of the above-mentioned methods consistently finds a near-optimal regularization parameter for all test problems and noise realizations in our numerical examples. In particular for rank-deficient problems the above-mentioned methods tend to produce solutions that significantly differ from the true solutions.

We can state the problem formally as follows. Given an ill-conditioned matrix \(A\) and vector \(b\) contaminated by noise, we solve the linear system
\be\label{eq:Ax=b} Ax=b.
\ee
Linear discrete ill-posed problems of the form \cref{eq:Ax=b} arise in a variety of settings, including the discretization of Fredholm integral equations of the first kind \cite{Kress2013,Hackbusch1995,Colton2013,Twomey1963,Hansen1992}, image deblurring problems \cite{SpecFiltBook,Dong2011,Yuan2007,Biemond1990,Oliveira2009}, machine learning algorithms \cite{Belkin2006,Vapnik1995,Smola2004,Cristianini2000,Williams2003} and more. The method of Tikhonov regularization replaces the original ill-posed problem \cref{eq:Ax=b} with a minimization problem
\be\label{eq:Pbx} \min_x\, ||Ax-b||^2+\lambda^2||Lx||^2,\ee
where \(||\cdot||\) is the \(\ell^2\)-norm, \(L\) is a regularization matrix and \(\lambda\) is a regularization parameter. The problem \cref{eq:Pbx} is said to be in standard form if \(L=I\), where \(I\) is the identity matrix, and in general form if \(L\neq I\) \cite{projRegGenForm,Hansen1989,Golub1999}. The development of an accurate and reliable method for determination of the regularization parameter \(\lambda\) is the main subject of this paper.

The present work is closely related to that of O'Leary \cite{PicInd2} and Taroudaki and O'Leary \cite{PicInd1}, who suggest a method for near-optimal estimation of the regularization parameter for standard-form Tikhonov regularization. The essence of this method is to determine the regularization parameter  by approximate minimization of the mean-square error (MSE)
\be\label{eq:MSEdefn}
\text{MSE}(\lambda) = ||x_{true} - x(\lambda)||^2,
\ee
where \(x_{true}\) is the least-squares solution of \(Ax=b_{true}\), \(b_{true}\) is the unperturbed data vector and \(x(\lambda)\) is the regularized solution of \cref{eq:Pbx}. Since the vector \(x_{true}\) is not known in practice, the authors of \cite{PicInd1} approximate the MSE \cref{eq:MSEdefn} by separating signal from noise using the coordinates of the perturbed data \(b\) with respect to the basis of the left singular vectors of \(A\), termed the Fourier coefficients of \(b\). Specifically, the authors of \cite{PicInd1} demonstrate the existence of an index they term the 'Picard parameter', which separates noise dominated coefficients of \(b\) from clean ones. Using the Picard parameter, the SVD expansion of the MSE is split into two parts - one containing the information about the unperturbed data and another containing the noise. The first part is replaced with its expected value, while the second part is rewritten in terms of the known Fourier coefficients of the data. The Picard parameter is either estimated manually in \cite{PicInd1}, representing the index where the plot of the Fourier coefficients of the data levels off, or, in the case where the noise is Gaussian, by the Lilliefors test, as we explain in \cref{sec:PictInd}.

Although the method developed in \cite{PicInd1,PicInd2} provides accurate results in a large percentage of cases and is shown to be competitive with standard methods, it holds two significant limitations. First, it does not allow the use of \(L\neq I\) which is necessary in many applications in order to incorporate various desirable properties in the solution \cite{Lchoice,DPC,projRegGenForm,Hansen1989}. In particular, \(L\) is often chosen to be the discrete approximation of a derivative operator to control various degrees of smoothness of the solution. The second, and arguably more important limitation is that the method is inaccurate for some noise realizations due to inaccurate estimation of the Picard parameter. This stems from the algorithm's reliance on applying statistical tests to noisy sequences, which often give inconsistent results.

In this paper we strive to overcome the above limitations. To handle the case \(L\neq I\), we replace the SVD of \(A\) used in \cite{PicInd1,PicInd2} with the generalized singular value decomposition (GSVD) of the pair \((A,L)\) \cite{GSVD,RankDeff}. It is significantly more difficult, however, to minimize the MSE \cref{eq:MSEdefn} in the GSVD basis \cite[sect. 2]{PicInd1}. For this reason, we replace the MSE with the predictive mean-square error (PMSE)
\be\label{eq:PMSEorig}
\text{PMSE}(\lambda) = ||A(x_{true}-x(\lambda))||^2,
\ee
where \(Ax(\lambda)\) is the data reconstructed from the regularized solution, termed the 'predicted data'. For simplicity, we assume that the unperturbed system is consistent, so that $b_{true}=Ax_{true}$, implying that the PMSE can be written as
\be\label{eq:PMSEdefn}
\text{PMSE}(\lambda) = ||b_{true}-Ax(\lambda)||^2.
\ee
The PMSE \cref{eq:PMSEdefn} has a simple expansion in terms of the GSVD basis, but it does not measure the error in the solution directly as the MSE. In principle, it is therefore possible that an algorithm successfully minimizing the PMSE will produce a suboptimal solution with high MSE in some problems. Nevertheless, the PMSE has approximately the same minimizer as the MSE in a variety of settings \cite{PMSEvsMSE1,PMSEvsMSE2} and the two minimizers were shown to coincide under certain assumptions \cite{PMSEvsMSE3}. To the best of our knowledge however, full characterization of the cases in which the minimizers of the MSE and PMSE are equal is unavailable. Therefore, we also provide a method for approximately minimizing the MSE as in \cite{PicInd1}, but for the more general case of \(L\neq I\). We then show that the expansion of the MSE in terms of the GSVD is numerically unstable and hence the accuracy of its approximation is limited. It is therefore advantageous to minimize the PMSE in cases where its minimizer is known to be close to that of the MSE, as it can be approximated better and hence leads to a better choice of \(\lambda\).

To determine \(\lambda\), we write the PMSE in terms of the GSVD of \((A,L)\) under very relaxed assumptions about the sizes and ranks of the matrices involved in \cref{eq:Pbx}. Next, we approximate \(\text{PMSE}(\lambda)\) by estimating the Picard parameter, splitting the GSVD expansion and modifying the noise-dependent terms as in \cite{PicInd1,PicInd2}. The regularization parameter \(\lambda\) is found by minimization of the resulting approximation of \(\text{PMSE}(\lambda)\). We term this procedure for determination of \(\lambda\) the Series Splitting (SS) method.

As an alternative to the SS method, approximate minimization of \cref{eq:PMSEdefn} can be performed using a more general two-step approach. Particularly, we can obtain an approximation \(\hat{b}\approx b_{true}\) by applying an accurate filter based on the Picard parameter to \(b\), substitute \(\hat b\) for \(b_{true}\) in \cref{eq:PMSEdefn} and minimize the resulting norm \(||\hat{b}-Ax(\lambda)||^2\). Our implementation of the filter employs the Picard parameter to remove the noise-dominated components of the data in the GSVD coordinate system. The advantage of this approach is its generalizability to other filters and regularization methods beside Tikhonov \cref{eq:Pbx}, as it requires only a data filter and an algorithm for calculating the regularized solution given a regularization parameter. We term this method the Data Filtering (DF) approach.

According to our numerical examples, the accuracy of the Picard parameter estimation algorithm \cite{PicInd1} is somewhat limited. We suggest that this can be significantly improved by a simple modification. Specifically, we propose to test the sequence of the Fourier coefficients in an order reverse to the one proposed in \cite{PicInd1} and at a higher confidence level. The performance of this algorithm can also be improved by providing upper and lower bounds on the Picard parameter, to limit the number of required tests. However, in spite of the improvement in accuracy, the algorithm remains prone to errors due to its reliance on noisy series of the Fourier coefficient of the data. To altogether avoid the dependence on noisy sequences, we proceed a step further to propose a new method that relies on averages of the squared moduli of the Fourier coefficients. We prove that the sequence of these averages decreases with increasing index of the Fourier coefficients, eventually converging to the value of the noise variance at the Picard parameter. In the resulting setting both the noise variance and the Picard parameter can be estimated reliably by detecting the levelling off of the above sequence of averages.

The SS method is closely related to the SURE and GCV methods, both of which approximately minimize \cref{eq:PMSEdefn}. In contrast to the SS and DF  methods, GCV and SURE do not split the sum in the GSVD expansion of \cref{eq:PMSEdefn}. Moreover, we show that both of them rely on replacing the whole sum in the expansion of \cref{eq:PMSEdefn} with its expected value, which results in an approximation less accurate than the one that could be achieved using the SS and DF methods. We provide a detailed comparison of the methods in a series of numerical examples.

The structure of this paper is as follows. In \cref{sec:prob} we formulate the problem of the Tikhonov regularization and solve it using the generalized singular value decomposition (GSVD). In \cref{sec:minPMSE}, we develop the SS and DF methods to approximately minimize the PMSE and the algorithms for estimation of the Picard parameter. In \cref{sec:others}, we discuss the SURE and GCV methods. Finally, in \cref{sec:numEx} we present the results of the numerical simulations.

\section{Formulation of the problem}\label{sec:prob}
We solve the linear ill-posed problem (\ref{eq:Ax=b}) by the Tikhonov regularization using a general regularization matrix \(L\). Throughout the paper, we make the following assumptions:
\begin{enumerate}
    \item\label{assump1} The problem \cref{eq:Pbx} has a unique solution for any \(\lambda\). This implies \(\mathcal{N}(A)\cap\mathcal{N}(L) = \{0\}\), where \(\mathcal{N}(\cdot)\) denotes the null-space of a matrix (see \cite[sect. 5.1.1]{RankDeff}).
    \item\label{assumpNL} The nullspace of \(L\) is spanned by smooth vectors. This assumption holds in most practical cases and is necessary to ensure proper filtering of the noise by the regularization, see \cite[sect. 8.1]{HansenInsights}, \cite[sect. 2.1.2]{RankDeff}, \cite[sect. 3]{Hansen1989}.
    \item\label{assump2} The data vector \(b\) is perturbed by an additive noise so \(b = b_{true} + n\), and the components of \(n\) are independent random variables taken from a normal distribution with zero mean and constant variance \(s^2\).
    \item The unperturbed system is consistent, so $b_{true} = Ax_{true}$.
    \item\label{assump3} The generalized singular values of \((A,L)\) decay to zero with no significant gap. The smallest generalized singular values cluster at (machine) zero. This property is common for discrete ill-posed problems - see \cite[sect. 2.1.2]{RankDeff}.
    \item\label{assump4} The problem satisfies the discrete Picard condition \cite{DPC}.
    \item\label{assump5} The minimizers of the MSE \cref{eq:MSEdefn} and the PMSE \cref{eq:PMSEdefn} are close to one another. This has been demonstrated in numerous numerical experiments such as \cite{PMSEvsMSE1,PMSEvsMSE2} and can be proved analytically for certain problems, as was done in \cite{PMSEvsMSE3} and \cite[Sect. 8.4]{splines}. However, to the best of our knowledge, no general characterization of cases in which the two minimizers are close is available. Therefore, for the completeness of this presentation, we provide a method for minimization of the MSE, in addition to the one minimizing the PMSE.
\end{enumerate}

The Tikhonov minimization problem \cref{eq:Pbx} is equivalent to the normal equation
\be\label{eq:problem} (A^*A+\lambda^2L^*L)x = A^*b, \ee
where \(^*\) denotes the conjugate transpose, thereby yielding the Tikhonov solution as
\be\label{eq:TikhSolnUseless} x(\lambda) = (A^*A+\lambda^2L^*L)^{-1}A^*b.\ee

We can express \cref{eq:TikhSolnUseless} in a more convenient form, using the GSVD \cite{GSVD} of the pair \((A,L)\).
To do so, let \(A\in \mathbb{C}^{m\times n}\), and \(L\in \mathbb{C}^{p\times n}\). Using these definitions, the GSVD of the matrices \(A\) and \(L\) is given by
\be\label{eq:GSVD}A = U\left(
                         \begin{array}{ccc}
                           I_A &  &  \\
                            & S_A &  \\
                            &  & O_A \\
                         \end{array}
                       \right)
Y^{-1}, \qquad L = V\left(
                      \begin{array}{ccc}
                        O_L &  &  \\
                         & S_L &  \\
                         &  & I_L \\
                      \end{array}
                    \right)
Y^{-1}, \ee
where
\begin{itemize}
    \item \(U\in \mathbb{C}^{m\times m}\), \(V\in \mathbb{C}^{p\times p}\) are unitary,
    \item \(Y \in \mathbb{C}^{n\times n}\) is invertible,
    \item \(S_A = \text{diag}\{\sigma_{r+1},\sigma_{r+2},...,\sigma_{r+q}\}\) and \(S_L = \text{diag}\{\mu_{r+1},\mu_{r+2},...,\mu_{r+q}\}\) are real diagonal matrices (note that \(\{\sigma_k\}\) are not the singular values of \(A\)),
    \item \(O_A \in \mathbb{C}^{(m-r-q)\times(n-r-q)}\), \(O_L \in \mathbb{C}^{(p+r-n)\times r}\) are zero matrices,
    \item \(I_A\) and \(I_L\) are \(r\times r\) and \((n-r-q)\times(n-r-q)\) identity matrices, respectively.
\end{itemize}
Note that the above zero and identity matrices can be empty. The values \(\{\sigma_k\}_{r+1}^{r+q}\) are arranged in decreasing order and \(\{\mu_k\}_{r+1}^{r+q}\) in increasing order so that
\be\label{eq:SingValArrange} 1>\sigma_{r+1}\geq\sigma_{r+2}\geq...\geq\sigma_{r+q}>0, \qquad 0<\mu_{r+1}\leq\mu_{r+2}\leq...\leq\mu_{r+q}<1.\ee
The pairs \((\sigma_k,\mu_k)\) satisfy the identity
\be\label{eq:SingValNormalize} \sigma_k^2+\mu_k^2 = 1\ \longleftrightarrow\ S_A^TS_A + S_L^TS_L = I_{q\times q}. \ee
The quantities \(\gamma_k = \sigma_k/\mu_k\) are the generalized singular values of the pair \((A,L)\). According to \cref{eq:SingValArrange} the sequence \(\{\gamma_k\}_{r+1}^{r+q}\) is arranged in decreasing order. Note, however, that the restriction of \(\sigma_k\) and \(\mu_k\) to \((0,1)\) does not apply to \(\gamma_k\).

To relate the parameters \(r\) and \(q\) to the ranks of \(A\) and \(L\), we observe that
\be A^*A = (Y^{-1})^*D_AY^{-1},\quad L^*L = (Y^{-1})^*D_LY^{-1},\ee
where \(D_A\) and \(D_L\) are diagonal matrices of the following form -
\be\label{eq:DaDl} D_A = \text{diag}\{\underbrace{1,...,1}_{r},\underbrace{\sigma_{r+1}^2,...,\sigma_{r+q}^2}_{q},\underbrace{0,...,0}_{n-r-q}\},\quad D_L = \text{diag}\{\underbrace{0,...,0}_{r},\underbrace{\mu_{r+1}^2,...,\mu_{r+q}^2}_{q},\underbrace{1,...,1}_{n-r-q}\}.\ee
Since \(D_A\) and \(D_L\) constitute diagonalizations of \(A^*A\) and \(L^*L\), it follows that the ranks of \(A^*A\) and \(L^*L\) are equal to the number of nonzero elements in \(D_A\) and \(D_L\), respectively, yielding the relations
\be\begin{dcases}\text{rank}(A) = r+q,\\
\text{rank}(L) = n-r,\end{dcases}\ \ee

For simplicity, we denote the columns of the matrices \(Y\) and \(U\) by \(\{y_k\}_{k=1}^{n}\) and \(\{u_k\}_{k=1}^{m}\), respectively. We also drop the indices from the column notations \(\{\cdot\}\) when referring to the entire column. The Fourier coefficients of the data and of the noise, with respect to the basis \(\{u_k\}\), are denoted by \(\beta_k= u_k^*b\) and \(\nu_k= u_k^*n\) respectively. Using these definitions and the decomposition \cref{eq:GSVD}, the Tikhonov solution can be written as
\be\label{eq:TikhSolnDisc}
x(\lambda) = \sum_{k=1}^r\beta_ky_k + \sum_{k=r+1}^{r+q} \frac{\gamma_k^2}{\gamma_k^2+\lambda^2}\frac{\beta_k}{\sigma_k}y_k.
\ee
If \(L\) is nonsingular, the generalized singular values \(\{\gamma_k\}\) are the regular singular values of \(AL^{-1}\) \cite[sect. 2.1.2]{RankDeff}, \cite{LCurveAnal}. In particular, if \(L=I\) the SVD of \(A\) is given by
\be\label{eq:SVDwithGSVD} A = U\left(
                                 \begin{array}{cc}
                                   S_AS_L^{-1} &   \\
                                     & O_A \\
                                 \end{array}
                               \right)V^*,\ee
where  \(r=0\), \(S_AS_L^{-1}=\text{diag}\{\gamma_1,\ldots,\gamma_q\}\) and the matrices \(O_A\), \(U\) and \(V\) are obtained from the GSVD \cref{eq:GSVD}. Furthermore, denoting the columns of \(V\) by \(\{v_k\}\) it is easy to show that \(y_j=\mu_j v_j\) for \(j\leq\text{rank}(A)\) and \(y_j=v_j\) for \(j>\text{rank}(A)\).
Thus, when \(L=I\) our expression for the Tikhonov solution \cref{eq:TikhSolnDisc} and the Fourier coefficients \(\beta_k\) and \(\nu_k\) coincide with the ones given in \cite{PicInd1,PicInd2}.

The factors $\gamma_k^2/(\gamma_k^2+\lambda^2)$ in \cref{eq:TikhSolnDisc} can be viewed as filters applied to the noisy data coefficients $\beta_k$, since they satisfy $\phi_{r+1} \approx 1$ and $\phi_{r+q}\approx 0$, and thereby dampen the coefficients for large $k$. While these coefficients correspond to the noise, the coefficients for small $k$ correspond to the true data and remain almost unchanged. In general, as discussed in \cite{PicInd1}, the coefficients can be replaced with more general filter factors $\phi_k(\lambda,\{\sigma_k\},\{\mu_k\})$, with a similar dampening effect. While we shall focus on the Tikhonov filters in this paper, all our subsequent derivations can be easily generalized to arbitrary filter factors such as those considered in \cite{PicInd1}.

\section{Estimation of the regularization parameter}\label{sec:minPMSE}
In this section we consider the problem of choosing a near-optimal value of the regularization parameter \(\lambda\) in \cref{eq:TikhSolnDisc}. The existence of such a value of \(\lambda\) is guaranteed by the discrete Picard condition \cite{DPC}, \cite[sect. 4.5]{RankDeff}, \cite[sect. 2.1]{PicInd1}, which requires the sequence \(\{|\beta_k-\nu_k|\}\) to decay faster than the generalized singular values \(\{\gamma_k\}\).
Since according to \cref{assump3}, we have \(\gamma_{j}\approx 0\) from some index \(j\) on, the discrete Picard condition implies the existence of an index \(k_0\leq j\) called the Picard parameter such that \(|\beta_k-\nu_k|<\gamma_k\approx 0\), or equivalently, \(\beta_k\approx \nu_k\) for all \(k\geq k_0\). This property of the Picard parameter \(k_0\) is used below to approximate the PMSE.

\subsection{The Series Splitting method}\label{sec:SS}
We assess the quality of \(x(\lambda)\) by measuring the distance \cref{eq:PMSEdefn}
between the unperturbed data \(b_{true}=b-n\) and the predicted data
\be\label{eq:Ax_l}
Ax(\lambda)  = \sum_{k=1}^r \beta_ku_k + \sum_{k=r+1}^{r+q} \frac{\gamma_k^2}{\gamma_k^2+\lambda^2}\beta_ku_k.\ee
We can rewrite \cref{eq:PMSEdefn} as
\beq\label{eq:PMSEnRewrite}
\text{PMSE}(\lambda) &=& ||b-n-Ax(\lambda)||^2 \nonumber\\
&=& ||n||^2+\rho(\lambda)-2C(\lambda),
\eeq
where
\be\label{eq:ResidualNorm}\rho(\lambda) = ||b-Ax(\lambda)||^2
= \sum_{k=r+1}^{r+q} \frac{\lambda^4}{(\gamma_k^2+\lambda^2)^2}|\beta_k|^2 + \sum_{k=r+q+1}^m|\beta_k|^2,\ee
is the squared residual norm,
\be\label{eq:TN} C(\lambda) = \Re[n^*(b-Ax(\lambda))] = \sum_{k=r+1}^{r+q}\frac{\lambda^2}{\gamma_k^2+\lambda^2}\Re(\beta_k\overline{\nu}_k) + \sum_{k=r+q+1}^m\Re(\beta_k\overline{\nu}_k),\ee
and \(\Re(\cdot)\) denotes the real part.

Noting that the term \(||n||^2\) in \cref{eq:PMSEnRewrite} is independent of \(\lambda\) and can therefore be neglected, we find that it is sufficient to minimize
\be\label{eq:gl} g(\lambda) = \rho(\lambda) -2C(\lambda).\ee
A direct evaluation of \(g(\lambda)\) is not possible since the function \(C(\lambda)\) depends on the unknown noise vector \(n\). Nonetheless, we can approximate \(C(\lambda)\) accurately using the Picard parameter \cite{PicInd1,PicInd2}. Recalling that the Picard parameter \(k_0\) is the smallest index for which  \(\beta_k\approx\nu_k\) is satisfied for all \(k\geq k_0\), we can split the sequence \(\{\beta_k\}\) into two parts - the first part \(\{\beta_k\}_{k=1}^{k_0-1}\), which contains the information about the unperturbed data and the second part \(\{\beta_k\}_{k=k_0}^{m}\), which contains the noise. Therefore, when \(k\geq k_0\) we can approximate the unknown term in \cref{eq:TN} by \(\Re(\beta_k\overline{\nu}_k)\approx|\beta_k|^2\). When \(k<k_0\) however, the coefficients \(\beta_k\) and \(\nu_k\) differ significantly, so we choose to approximate the term \(\Re(\beta_k\overline{\nu}_k)\) in \cref{eq:TN} by replacing it with its expected value. Denoting the expected value by \(\mathcal{E}(\cdot)\), and noting that \(\beta_k=(\beta_k-\nu_k)+\nu_k\) where \(\beta_k-\nu_k = u_k^*b_{true}\) is not random, we can deduce using \cref{assump2} that
\be\label{eq:Expect1}
\begin{aligned} \mathcal{E}(\Re(\beta_k\overline{\nu}_k)) &= \Re\left(\mathcal{E}(\beta_k\overline{\nu_k})\right) \\
&= \Re\left(\mathcal{E}[(\beta_k-\nu_k)\overline{\nu}_k]\right) + \Re\left(\mathcal{E}(|\nu_k|^2)\right) \\
&= \Re\left((\beta_k-\nu_k)\underbrace{\overline{\mathcal{E}(\nu_k)}}_{=0}\right) + s^2\\ &= s^2,\end{aligned}
\ee
see \cite[sect. 6.6]{SpecFiltBook}, \cite{PicInd1,PicInd2}. Therefore, for given \(k_0\) and \(s^2\) we can approximate \(C(\lambda)\) by splitting the series \cref{eq:TN} similarly to \cite{PicInd1} to obtain
\be\label{eq:TnApprox} C(\lambda) \approx s^2\sum_{k=r+1}^{k_0-1}\frac{\lambda^2}{\gamma_k^2+\lambda^2} +\sum_{k=k_0}^{r+q}\frac{\lambda^2}{\gamma_k^2+\lambda^2}|\beta_k|^2 +\sum_{k=r+q+1}^{m}|\beta_k|^2.\ee
The regularization parameter is then found by minimizing \cref{eq:gl} using the approximation of \(C(\lambda)\) given by \cref{eq:TnApprox}.

We can show that \(k_0\) is limited to the interval \([r+1, r+q]\) and therefore it is always possible to split $C(\lambda)$ as in \cref{eq:TnApprox}. To justify the lower bound, we note that the nullspace of \(L\) is spanned by the vectors \(\{y_k\}_{k=1}^r\), as they constitute a set of \(r=\dim\mathcal{N}(L)\) linearly independent vectors satisfying \(Ly_k=0\). Since the vectors \(\{y_k\}_{k=1}^r\) are smooth (by \cref{assumpNL}) and \(A\) has a typical smoothing effect \cite[p. 21]{RankDeff}, \(\{u_k\}_{k=1}^r\) are also smooth and satisfy \(Ay_k = u_k\) for \(k\leq r\). Therefore, the smooth vector \(b_{true}\) is well-represented by vectors \(\{u_k\}_{k=1}^r\), while the non-smooth noise vector \(n\) is represented mostly by \(u_k\) with \(k>r\). Thus, we have \(\nu_k\ll \beta_k\) for \(k\leq r\), implying that \(\beta_k-\nu_k\not\approx 0\) for \(k\leq r\) and so, \(k_0\geq r+1\). To justify the upper-bound, we note that \(\gamma_{r+q}\) is the last generalized singular value of \(A\) which is numerically nonzero and, by \cref{assump3}, \(\gamma_{r+q}\approx \epsilon\) where \(\epsilon\) is the machine zero. Thus, by the discrete Picard condition, \(\beta_{r+q}-\nu_{r+q}\approx 0\) and we can conclude that \(k_0\leq r+q\).
If \(L=I\) as in \cite{PicInd1,PicInd2}, we have \(r=0\) and therefore only the upper bound \(k_0\leq\text{rank}(A)\) is nontrivial.
The SS algorithm is summarized in \cref{alg:SSAlg}.

\subsection{Approximate minimization of the MSE} As discussed above, our approach approximately minimizes the PMSE \cref{eq:PMSEdefn}, in contrast to the approach of \cite{PicInd1,PicInd2} which minimizes the MSE \cref{eq:MSEdefn} assuming \(L=I\). To illustrate the difference between the two approaches and to give an alternative method of solution for problems whose PMSE and MSE minimizers may not coincide, we repeat the above derivation for the MSE and approximate it in the general case, \(L\neq I\).

We begin with the observation that \(x_{true}\), the least-squares solution to \(Ax=b_{true}\) can be expressed in the GSVD basis as
\be\label{eq:TrueSolnGSVD}
x_{true} = \sum_{k=1}^r(\beta_k-\nu_k)y_k + \sum_{k=r+1}^{r+q} \frac{\beta_k-\nu_k}{\sigma_k} y_k.
\ee
We can rewrite \cref{eq:TrueSolnGSVD} as
\be\label{eq:TrueSolnGSVDNotation}
x_{true} = x(0) - \widehat x,
\ee
where \(x(0)=x_{LS}\) is the least-squares solution of the perturbed problem \cref{eq:Ax=b} obtained by substituting \(\lambda=0\) in \cref{eq:TikhSolnDisc} and
\[
\widehat x=\sum_{k=1}^r\nu_k y_k + \sum_{k=r+1}^{r+q}\frac{\nu_k}{\sigma_k}y_k,
\]
is the least-squares solution for the pure noise problem \(Ax=n\). The MSE can then be expanded as
\be\label{eq:MSEexpand}
\text{MSE}(\lambda) = \widetilde\rho(\lambda) + ||\widehat x||^2 - 2D(\lambda),
\ee
where
\be\label{eq:residualNormRecSpace}
\widetilde\rho(\lambda)= ||x(0)-x(\lambda)||^2 = \sum_{j=r+1}^{r+q}\sum_{k=r+1}^{r+q}\frac{y_j^*y_k}{\sigma_j\sigma_k}\frac{\lambda^4}{(\gamma_j^2+\lambda^2)(\gamma_k^2+\lambda^2)}\beta_k\overline{\beta}_j,
\ee
and
\be\label{eq:DN}
\begin{aligned} D(\lambda) = \Re\left(\widehat x^*\left(x(0)-x(\lambda)\right)\right) = &\sum_{j=1}^r\sum_{k=r+1}^{r+q}\frac{y_j^*y_k}{\sigma_k} \frac{\lambda^2}{\gamma_k^2+\lambda^2}\Re(\beta_k\overline{\nu}_j)\\ &+ \sum_{j=r+1}^{r+q}\sum_{k=r+1}^{r+q} \frac{y_j^*y_k}{\sigma_j\sigma_k} \frac{\lambda^2}{\gamma_k^2+\lambda^2}\Re(\beta_k\overline{\nu}_j).\end{aligned}\ee
The first term in \cref{eq:MSEexpand}, \(\widetilde\rho(\lambda)\), can be readily evaluated while the second term, \(||\widehat x||^2\), can be dropped entirely as it does not depend on \(\lambda\). The third term, \(D(\lambda)\) cannot be evaluated as it depends on the coefficients \(\nu_k\) of the unknown noise vector \(n\) as shown in \cref{eq:DN} and must therefore be approximated.

To approximate \(D(\lambda)\) we first rewrite the sums in \cref{eq:DN} as
\be\label{eq:DNsplit}
\begin{aligned} D(\lambda) =  &\sum_{j=1}^r\sum_{k=r+1}^{r+q}\frac{y_j^*y_k}{\sigma_k} \frac{\lambda^2}{\gamma_k^2+\lambda^2}\Re(\beta_k\overline{\nu}_j) + \sum_{k=r+1}^{k_0-1}\frac{||y_k||^2}{\sigma_k^2} \frac{\lambda^2}{\gamma_k^2+\lambda^2}\Re(\beta_k\overline{\nu}_k)\\ &+ \left( \sum_{j=r+1}^{k_0-1}\sum_{\substack{k=r+1\\k\neq j}}^{k_0-1} +\sum_{j=k_0}^{r+q}\sum_{k=k_0}^{r+q}\right)\frac{y_j^*y_k}{\sigma_j\sigma_k} \frac{\lambda^2}{\gamma_k^2+\lambda^2}\Re(\beta_k\overline{\nu}_j).\end{aligned}\ee
When \(j\geq k_0\), we can approximate the coefficients of the noise, \(\nu_j\), by the coefficients of the data, \(\beta_j\), so that \(\Re\left(\beta_k\overline{\nu}_j\right)\approx\beta_k\overline{\beta}_j\). However, when \(j<k_0\) we approximate terms involving \(\overline{\nu}_j\) in \cref{eq:DNsplit} by replacing them with their expected value as we do for \(C(\lambda)\) in \cref{eq:TnApprox}. The main difference is that \cref{eq:TnApprox}, in addition to the case \(j=k\) for which \(\mathcal{E}(\Re(\beta_k\overline{\nu}_k)) = s^2\) as in \cref{eq:Expect1}, contains cross terms with \(j\neq k\). However these terms can be neglected since \cref{assump2} implies
\[
\mathcal{E}(\Re(\beta_k\overline{\nu}_j)) = \Re\left(\mathcal{E}(\nu_k\overline{\nu_j})\right)=\Re\left(\mathcal{E}(\nu_k)\mathcal{E}(\overline{\nu_j})\right)=0.
\]
Consequently, we can drop the first and third sum in \cref{eq:DNsplit}, approximate \(D(\lambda)\) as
\be\label{eq:DNapprox}
D(\lambda)\approx \widetilde D(\lambda) = s^2\sum_{k=r+1}^{k_0-1}\frac{||y_k||^2}{\sigma_k^2} \frac{\lambda^2}{\gamma_k^2+\lambda^2} + \sum_{j=k_0}^{r+q}\sum_{k=k_0}^{r+q} \frac{y_j^*y_k}{\sigma_j\sigma_k} \frac{\lambda^2}{\gamma_k^2+\lambda^2}\beta_k\overline{\beta}_j,
\ee
and estimate the minimum of the MSE in \cref{eq:MSEexpand} by minimizing
\be\label{eq:GlSS}
\widetilde g(\lambda) = \widetilde\rho(\lambda) - 2\widetilde D(\lambda).
\ee
The problem with this approach is that \cref{eq:residualNormRecSpace} and \cref{eq:DNapprox} are numerically unstable due to the division by \(\sigma_k\sigma_j\) for \(k,j\geq k_0\). Specifically, due to the ill-conditioning of \(A\), the values \(\{\sigma_j\}\) decay quickly to zero\footnote{More precisely, \cref{assump4} requires \(\gamma_{r+q}\approx\epsilon\), but since \(\gamma_{r+q}=\sigma_{r+q}/\mu_{r+q}\) and \(\mu_{r+q}=\sqrt{1-\sigma_{r+q}^2}\) by \cref{eq:SingValNormalize}, we have \(\sigma_{r+q}\approx \epsilon\sqrt{1-\epsilon^2}\approx \epsilon\).} so that \(\sigma_{r+q}\approx\epsilon\). Therefore, terms that include division by \(\sigma_k\sigma_j\) for large \(j\) and \(k\) completely dominate the value of \(\widetilde g(\lambda)\) and, due to finite machine precision, eliminate the contribution of the terms with small \(j\) and \(k\). This contrasts with the fact that for a desirable choice of \(\lambda\), the solution \(x(\lambda)\) should be smooth and its error should therefore depend significantly on terms with \(j,k<k_0\). To circumvent this problem, we drop the terms with \(j,k\geq k_0\) in the sums \cref{eq:residualNormRecSpace} and \cref{eq:DNapprox}, similarly to \cite{PicInd1}, so that
\be\label{eq:resNormRecSpaceTrunc}
\widetilde\rho(\lambda) \approx \widehat\rho(\lambda) = \sum_{j=r+1}^{k_0-1}\sum_{k=r+1}^{k_0-1}\frac{y_j^*y_k}{\sigma_j\sigma_k}\frac{\lambda^4}
{(\gamma_j^2+\lambda^2)(\gamma_k^2+\lambda^2)}\beta_k\overline{\beta}_j,
\ee
\be\label{eq:DNapproxTunc}
\widetilde D(\lambda) \approx \widehat D(\lambda) = s^2\sum_{k=1}^{k_0-1}\frac{||y_k||^2}{\sigma_k^2} \frac{\lambda^2}{\gamma_k^2+\lambda^2},
\ee
and minimize \(\widehat g(\lambda) = \widehat\rho(\lambda)-2\widehat D(\lambda)\), thereby retaining only terms in which \(\sigma_k\) is relatively large. Note that if \(L=I\), the above method of minimizing \(\widehat g(\lambda)\) becomes identical to Algorithm 1 in \cite{PicInd1} (implemented with the Tikhonov filter factors).

In spite of the resulting numerical stability of the above scheme upon dropping of terms with \(j,k\geq k_0\), we lose information that might have improved the accuracy of the approximation of the MSE if there was no instability. In contrast, no such division by \(\sigma_j\sigma_k\) is necessary for minimization of \(\text{PMSE}(\lambda)\) \cref{eq:PMSEdefn}, in which case we can retain all the terms upon the series splitting in \cref{eq:TnApprox}.
Thus, we expect our approximation \cref{eq:gl} of norm (1.4) to yield a better estimate of \(\lambda\) compared to the approximation minimizing \(\widehat g(\lambda)\), where there is no discrepancy between the minimal of MSE and PMSE. This situation is quite common, as it is shown in \cite[Sect. 8.4]{splines} and \cite{PMSEvsMSE3} that the PMSE in \cref{eq:PMSEdefn} and the MSE \cref{eq:MSEdefn} have approximately the same minimizer in a variety of settings, and numerical results supporting this are available in e.g. \cite{PMSEvsMSE1,PMSEvsMSE2}. We therefore focus below on minimizing the PMSE and not the MSE, and demonstrate in the numerical examples of \cref{sec:numEx} that this approach yields superior results.

\subsection{Estimating the Picard parameter and the variance of the noise}\label{sec:PictInd}
We begin with a brief discussion of the methods for estimation of the Picard parameter \(k_0\), suggested in \cite{PicInd1}. The Picard parameter can be graphically deduced from the plot of the sequence \(\{\beta_k\}\) versus \(k\). Specifically, due to the discrete Picard condition, the plot of \(\{\beta_k\}\) is expected to decay on average with increasing index and to level-off at the Picard parameter. This levelling-off can be found manually from the plot, see \cite[sect. 2.2]{PicInd1}. The drawback of this approach is that due to the significant variance of plot \(\{\beta_k\}\) the point at which the plot levels-off cannot be unambiguously determined. In order to reduce this ambiguity, as well as to automate the method, it is suggested in \cite[sect. 2.3]{PicInd1} to use the Lilliefors test for normality on subsequences of \(\{\beta_k\}\). Specifically, this method sets \(k_0\) to the smallest index for which the sequence \(\{\beta_k\}_{k=k_0}^{m}\) is dominated by Gaussian noise. By applying the Lilliefors test at 95\% confidence to the sequences \(\{\beta_k\}_{k=j}^{m}\) for \(j=m-3,m-2,...,1\), the  Picard parameter \(k_0\) is chosen to be the smallest index after which the test fails 10 consecutive times. If the test fails immediately at \(j=m-3\), it is proposed to set \(k_0=m+1\) and \(s^2=0\), which signifies that the data is noiseless. Alternative tests that assume distributions different from the Gaussian distribution can be utilized in a similar way \cite{PicInd1}. Once \(k_0\) is found, the variance \(s^2\) can be estimated as the sample variance of the sequence \(\{\beta_k\}_{k=k_0}^m\) using the expression
\be\label{eq:estS2}
s^2\approx \frac{1}{m-k_0+1}\sum_{k=k_0}^m\left|\beta_k\right|^2,
\ee
where we note that, according to \cref{assump2}, the noise terms \(\{\nu_k\}\) have a zero mean, and therefore the mean of the sequence \(\{\beta_k\}_{k=k_0}^m\approx\{\nu_k\}_{k=k_0}^m\) is negligible.

One can improve the accuracy of the method described above by initializing the estimate of \(k_0\) to its lower bound \(r+1\) and applying the Lilliefors test to the sequences \(\{\beta_k\}_{k=j}^{m}\) for \(j=r+1,r+2,...,r+q\) at a 99.9\% confidence level. The value of \(k_0\) is set to the smallest index \(j\) for which the Lilliefors test indicates that the sequence is normally distributed. If the test fails for \(j=r+q\), we set \(k_0=r+q\) and \(s^2=0\). Once \(k_0\) is estimated, we can find the variance using \cref{eq:estS2}. This modified algorithm is summarized in \cref{alg:PicIndLilAlg}.

\begin{figure}
\centering
\includegraphics[width=.6\textwidth]{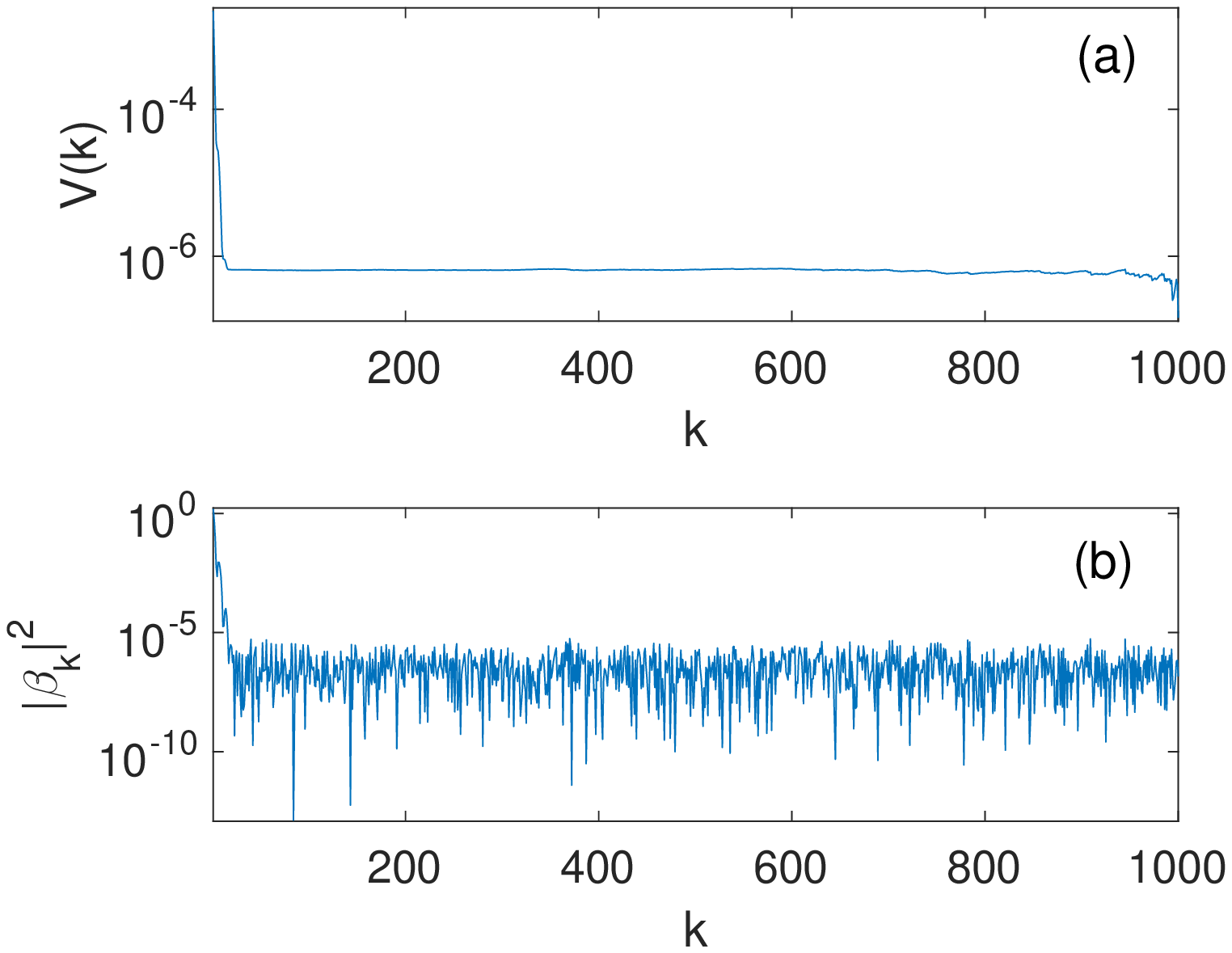} 
\caption{Comparison of (a) \(\log V(k)\) and (b) \(\log |\beta_k|^2\) for the test problem \texttt{heat} from \cite{RegTools}, corrupted with white Gaussian noise with \(s^2/\max\{|b_{true ; k}|^2\} = 10^{-4}\). Here, we use \(L=I\) as in \cite{PicInd1}. The non-negligible fluctuations of \(|\beta_k|^2\) in (b) and in contrast, almost flat \(V(k)\) in (a) are clearly seen.}\label{fig:VxVSWx}
\end{figure}

The dependence of the estimation of the Picard parameter upon statistical tests can be avoided using the following new method. This method is based on an averaging of the Fourier coefficients, which reduces the variance of the sequence \(\{|\beta_k|^2\}\) (see \cref{fig:VxVSWx}(a)), enabling more reliable automatic detection of the levelling-off of these coefficients. We note first that the sequence \(\{|\beta_k|^2\}\) decreases on average until it levels-off at \(k_0\) and oscillates about \(s^2\) with a non-negligible variance. To show this, we observe that
\be\label{eq:ExpOfbSq}\begin{aligned} \mathcal{E}(|\beta_k|^2) &= \mathcal{E}\left(|\beta_k-\nu_k|^2\right) - \underbrace{\mathcal{E}\left(|\nu_k|^2\right)}_{=s^2} + 2\underbrace{\mathcal{E}\left(\Re\left(\beta_k\overline{\nu_k}\right)\right)}_{=s^2} \\ &= |\beta_k-\nu_k|^2 + s^2.\end{aligned}\ee
Therefore, due to the discrete Picard condition, the expected value \(\mathcal{E}(|\beta_k|^2)\) must decrease on average with increasing \(k\) and become constant at \(s^2\) for \(k\geq k_0\). While the actual curve of \(|\beta_k|^2\) deviates from its expected value, these deviations are random and thus the curve of \(|\beta_k|^2\) oscillates about its expected value. This implies that the general trend of \(\mathcal{E}(|\beta_k|^2)\) to decrease on average and to level-off also applies to \(|\beta_k|^2\), which, in its turn, justifies the graphical method of \cite{PicInd1}. However, instead of flattening at \(k_0\) as \(\mathcal{E}(|\beta_k|^2)\), the curve of \(|\beta_k|^2\) oscillates about \(s^2\) for \(k\geq k_0\) with a non-negligible variance given by
\be\label{eq:VarOfBeta}
\text{Var}\left(|\beta_k|^2\right) \approx \mu_4 - s^4,\quad \text{for } k\geq k_0,
\ee
where \(\mu_4=\mathcal{E}(|\beta_k|^4)=\mathcal{E}(|\nu_k|^4)\) is the fourth moment of the noise distribution. For the derivation of \cref{eq:VarOfBeta} we again use \(\beta_k\approx\nu_k\) for \(k\geq k_0\). In particular, for Gaussian noise we have \(\mu_4=3s^4\) and therefore \(\text{Var}(|\beta_k|^2)=2s^4\). Thus, due to the significant variance of the sequence \(\{|\beta_k|^2\}_{k=k_0}^m\), any estimation of \(k_0\) from the levelling-off of \(\{|\beta_k|^2\}\) is prone to error. Therefore, instead of \(\{|\beta_k|^2\}\), we consider the sequence of averages \(\{V(k)\}\) given by
\be\label{eq:VxDefn}
V(k) = \frac{1}{m-k+1}\sum_{j=k}^m|\beta_j|^2.
\ee
Since the sequence \(\{|\beta_k|^2\}\) decays on average, so does the sequence \(\{V(k)\}\). To demonstrate this, we note that
\be\label{eq:ExpVk}
\mathcal{E}(V(k)) = \frac{1}{m-k+1}\sum_{j=k}^{k_0}|\beta_j-\nu_j|^2 + s^2,
\ee
and since \(\{|\beta_j-\nu_j|\}\) decays on average, so does \(V(k)\), as follows from the inequality
\be
\frac{1}{m-k_1+1}\sum_{j=k_1}^{k_0}|\beta_j-\nu_j|^2 > \frac{1}{m-k_2+1}\sum_{j=k_2}^{k_0}|\beta_j-\nu_j|^2, \quad \text{for } k_1<k_2<k_0.
\ee
In addition, \cref{eq:ExpVk} implies that \(V(k)\approx s^2\) for \(k\geq k_0\) and therefore the sequence \(\{V(k)\}\) levels-off at \(k_0\) and oscillates about \(s^2\), similar to \(\{|\beta_k|^2\}\). However, the variance of \(V(k)\) for \(k\geq k_0\) is significantly smaller compared to that of \(|\beta_k|^2\), making its curve significantly flatter and more suitable for the estimation of \(k_0\) and \(s^2\). Specifically, the variance of \(V(k)\) for \(k\geq k_0\) is given by
\be\label{eq:varOfVar}
\text{Var}(V(k)) = \frac{\mu_4-s^4}{m-k+1} = \frac{1}{m-k+1}\text{Var}(|\beta_k|^2),\quad \text{for } k\geq k_0,
\ee
and for the Gaussian noise it is \[ \text{Var}(V(k)) = \frac{2s^4}{m-k+1}.\]
Since the variance decays as \(1/(m-k+1)\)  the curve of \(V(k)\) remains practically flat for a wide range of indices \(k\).

In \cref{fig:VxVSWx} we illustrate the difference between the sequences \(\{|\beta_k|^2\}\) and \(\{V(k)\}\) for estimation of the Picard parameter. Even though both sequences decrease on average until they level-off and oscillate about \(s^2\), the plot of \(\{V(k)\}\) shown in \cref{fig:VxVSWx}(a) remains almost constant for \(k\geq k_0\), whereas \(\{|\beta_k|^2\}\) shown in \cref{fig:VxVSWx}(b) oscillates with a significant variance. Due to these oscillations the exact point at which the plot of  \(\{|\beta_k|^2\}\) levels-off cannot be unambiguously determined. In contrast, \(V(k)\) in \cref{fig:VxVSWx}(a) flattens almost completely and therefore the point at which it levels-off is easily found to be \(k_0\approx17\). Note that the nonzero variance of \(V(k)\) becomes significant at about \(k\approx 940\), where \(\text{Var}\left(V(k)\right)\approx s^4/40\).

To estimate \(k_0\) from the flatness of \(V(k)\) we suggest the following simple rule. We set \(k_0\) to the smallest index for which the relative change in \(V(k)\) is small enough, so that
\be\label{eq:PicIndCond1}
\frac{|V(k+h)-V(k)|}{V(k)} < \varepsilon,
\ee
for some bound \(\varepsilon\), where we require \(k\leq m-h\) in \cref{eq:PicIndCond1}. If \(V(k)\) does not satisfy \cref{eq:PicIndCond1} for any \(r+1\leq k_0\leq \min(r+q,m-h)\) we set \(k_0=r+q\) and \(s^2=0\). The upper bound of \(k_0\) arises because it should satisfy \(k_0\leq r+q\), as discussed in \cref{sec:minPMSE}, and \(k_0\leq m-h\) due to \cref{eq:PicIndCond1}. Once \(k_0\) is found using \cref{eq:PicIndCond1}, we estimate the variance as \(s^2\approx V(k_0)\) which follows from \cref{eq:estS2}.

Condition \cref{eq:PicIndCond1} depends on two free parameters, \(h\) and \(\varepsilon\). Step size \(h\) must be large enough to ensure that the flattening of \(V(k)\) is not due to random oscillations. However, it also needs to be small enough, \(h < r+q-k_0\), so as not to exclude the flat part of \(V(k)\) from consideration. Similarly, the value of \(\varepsilon\) has to be small enough to detect the levelling off of \(V(k)\) but large enough to account for its small but nonzero variance. In practice however, we found the criterion \cref{eq:PicIndCond1} to be robust to the choice of \(\varepsilon\) and \(h\).
We summarize the above procedure for estimation of \(k_0\) and \(s^2\) in \Cref{alg:PicIndAlg}. Note that the estimate of \(k_0\) from \Cref{alg:PicIndAlg} is used as input into \cref{alg:SSAlg} of \cref{sec:SS}.

\begin{algorithm}
\caption{The SS method}\label{alg:SSAlg}
\begin{algorithmic}
\INPUT{\(A,L,b,\varepsilon,h\)}
\OUTPUT{\(x_{SS}\)}
\State \([m, n]\gets \text{size}(A)\) \Comment{\(\text{size}\) = [no. of rows, no. of columns]}
\State \(r\gets n-\text{rank}(L)\)
\State \(q\gets \text{rank}(A)+\text{rank}(L)-n\)
\State \([\{\sigma_k\},\{\mu_k\},\{y_k\},\{u_k\}] \gets \text{GSVD}(A,L\)) \Comment{Perform GSVD of the pair \((A,L)\)}
\State \(\{\gamma_k\}\gets\{\sigma_k/\mu_k\}\) \Comment{Define the generalized singular values}
\State \(\{\beta_k\}\gets\{u_k^*b\}\) \Comment{Define the Fourier coefficients of \(b\) with respect to \(\{u_k\}\)}
\State \(x(\lambda) = \sum_{k=1}^r\beta_ky_k + \sum_{k=r+1}^{r+q} \frac{\gamma_k^2}{\gamma_k^2+\lambda^2}\frac{\beta_k}{\sigma_k}y_k\) \Comment{Define the Tikhonov solution as in \cref{eq:TikhSolnDisc}}
\State \(\rho(\lambda) = ||b-Ax(\lambda)||^2\) \Comment{Define \(\rho(\lambda)\) as in \cref{eq:ResidualNorm}}
\LineComment{Get \(k_0\) and \(s^2\) from \cref{alg:PicIndAlg} or alternative algorithms such as \cref{alg:PicIndLilAlg}:}
\State \([k_0,s^2] \gets \text{\cref{alg:PicIndAlg}}(\{\beta_k\},r,q,\varepsilon,h)\)
\State \(\widehat C(\lambda) = s^2\sum_{k=r+1}^{k_0-1}\frac{\lambda^2}{\gamma_k^2+\lambda^2} +\sum_{k=k_0}^{r+q}\frac{\lambda^2}{\gamma_k^2+\lambda^2}|\beta_k|^2 +\sum_{k=r+q+1}^{m}|\beta_k|^2\) \Comment{Define \(\widehat C(\lambda)\) as in \cref{eq:TN}}
\State \(g(\lambda) = \rho(\lambda) -2\widehat C(\lambda)\) \Comment{Define the function to be minimized as in \cref{eq:gl}}
\State \(\lambda_{found}\gets\arg\min_{\lambda} g(\lambda)\) \Comment{Find \(\lambda\) minimizing \(g(\lambda)\)}
\State \(x_{SS}\gets x(\lambda_{found})\)
\end{algorithmic}
\end{algorithm}

\begin{algorithm}
\caption{Estimating the Picard parameter using the Lilliefors test}\label{alg:PicIndLilAlg}
\begin{algorithmic}
\INPUT{\(\{\beta_k\},r,q\)}
\OUTPUT{\(k_0,s^2\)}
\State \(k_0\gets r+1\) \Comment{Initialize the Picard parameter \(k_0\) to its lower bound}
\While{\(\texttt{lilliefors}\left(\{\beta_k\}_{k=k_0}^m,0.999\right)\neq 0\) \textbf{and}\ \(k_0<r+q+1\)}
\LineComment{\texttt{lilliefors}(\(\{\cdot\}\),.999) returns 1 if the Lilliefors test rejects the null hypothesis at \(99.9\%\)}\par \LineComment{  confidence level and 0 otherwise}
\State \(k_0\gets k_0+1\) 
\EndWhile
\If{\(k_0=r+q+1\)} \Comment{Check whether \(k_0\) exceeds the upper bound}
\State \(s^2\gets 0\) \Comment{If yes, data is noiseless}
\State \(k_0\gets r+q\)
\Else
\State \(s^2\gets \left(\sum_{k=k_0}^{m}|\beta_k|^2\right)/(m-k_0+1)\) \Comment{Otherwise estimate the variance \(s^2\) as in \cref{eq:estS2}}
\EndIf
\end{algorithmic}
\end{algorithm}

\begin{algorithm}
\caption{Estimating the Picard parameter using the sequence \(\{V(k)\}_{k=r+1}^{k=r+q}\)}\label{alg:PicIndAlg}
\begin{algorithmic}
\INPUT{\(\{\beta_k\},r,q,\varepsilon,h\)}
\OUTPUT{\(k_0,s^2\)}
\State \(k_0\gets r+1\) \Comment{Initialize the Picard parameter \(k_0\) to its lower bound}
\State \(V(k)= \left(\sum_{j=k}^{m}|\beta_j|^2\right)/(m-k+1)\) \Comment{Define \(V(k)\) as in \cref{eq:VxDefn}}
\While{\(|V(k_0+h)-V(k_0)|/V(k_0+1)\geq\varepsilon\)\ \textbf{and}\ \(k_0<r+q-h+1\)}
\State \(k_0\gets k_0+1\) \Comment{Increase \(k_0\) as long as the condition \cref{eq:PicIndCond1} is not satisfied and \(k_0\) is below}\par \Comment{the upper bound}
\EndWhile
\If{\(k_0=r+q+1\)} \Comment{Check whether \(k_0\) exceeds the upper bound}
\State \(s^2\gets 0\) \Comment{If yes, data is noiseless}
\State \(k_0\gets r+q\)
\Else
\State \(s^2\gets V(k_0)\) \Comment{Otherwise estimate the variance \(s^2\) as \(V(k_0)\)}
\EndIf
\end{algorithmic}
\end{algorithm}

\subsection{The Data Filtering method}
In this section we describe the Data Filtering (DF) method for minimization of \cref{eq:PMSEdefn}, which generalizes the SS method. In the DF method we minimize the norm
\be\label{eq:DFdist}
\hat f(\lambda)=||\hat b - Ax(\lambda)||^2,
\ee
where \(\hat b\) is the filtered perturbed data. The Picard parameter is used to directly approximate the true data \(b-n\approx \hat b\), instead of approximating the noise-dependent terms in \cref{eq:PMSEdefn}, as done in \cref{sec:SS}. We assume that the sequence \(\{\beta_k\}_{k=1}^{k_0-1}\) is dominated by the signal and can be regarded as the true data in basis \(\{u_k\}\), while the sequence \(\{\beta_k\}_{k=k_0}^{m}\) is dominated by noise. To approximate \(b-n\) we drop the noise-dominated terms from the expansion of \(b\) in terms of the basis \(\{u_k\}\) to obtain the approximation
\be\label{eq:bTApprox}
\hat b = \sum_{k=1}^{k_0-1}\beta_ku_k.
\ee
The norm \cref{eq:DFdist} can then be written as
\be\label{eq:DFdist2}
\hat f(\lambda)= \sum_{j=r+1}^{k_0-1}\frac{\lambda^4}{(\gamma_j^2+\lambda^2)^2}\beta_j^2 + \sum_{j=k_0}^{r+q}\frac{\gamma_j^4}{(\gamma_j^2+\lambda^2)^2}\beta_j^2.
\ee
The DF method is summarized in \cref{alg:DFAAlg}. Note that the DF method can be generalized to use various data filters to obtain \(\hat b\) and various regularization methods to obtain \(x(\lambda)\). Therefore, minimization of norm \cref{eq:DFdist} represents a new general approach for the estimation of \(\lambda\) and will be addressed in more detail in a forthcoming work.

\begin{algorithm}
\caption{The DF method}\label{alg:DFAAlg}
\begin{algorithmic}
\INPUT{\(A,L,b,\varepsilon,h\)}
\OUTPUT{\(x_{DFA}\)}
\State \([m, n]\gets \text{size}(A)\) \Comment{\(\text{size}\) = [no. of rows, no. of columns]}
\State \(r\gets n-\text{rank}(L)\)
\State \(q\gets \text{rank}(A)+\text{rank}(L)-n\)
\State \([\{\sigma_k\},\{\mu_k\},\{y_k\},\{u_k\}] \gets \text{GSVD}(A,L\)) \Comment{Perform GSVD of (A,L)}
\State \(\{\gamma_k\}\gets\{\sigma_k/\mu_k\}\) \Comment{Define the generalized singular values}
\State \(\{\beta_k\}\gets\{u_k^*b\}\) \Comment{Define the Fourier coefficients of \(b\) with respect to \(\{u_k\}\)}
\LineComment{Get \(k_0\) from \cref{alg:PicIndAlg} or from alternative algorithms such as \cref{alg:PicIndLilAlg} or the algorithm of \cite{PicInd1}.}
\State \(k_0 \gets \text{\cref{alg:PicIndAlg}}(\{\beta_k\},r,q,\varepsilon,h)\)
\State \(\hat b \gets \sum_{k=1}^{k_0-1}\beta_ku_k\) \Comment{Get filtered data as in \cref{eq:bTApprox}}
\State \(x(\lambda) = \sum_{k=1}^r\beta_ky_k + \sum_{k=r+1}^{r+q} \frac{\gamma_k^2}{\gamma_k^2+\lambda^2}\frac{\beta_k}{\sigma_k}y_k\) \Comment{Define Tikhonov solution as in \cref{eq:TikhSolnDisc}}
\State \(\hat f(\lambda) = ||\hat b-Ax(\lambda)||^2\) \Comment{Define the function to be minimized}
\State \(\lambda_{found}\gets\arg\min_{\lambda} \hat f(\lambda)\) \Comment{Find \(\lambda\) minimizing \(\hat f(\lambda)\)}
\State \(x_{DFA}\gets x(\lambda_{found})\)
\end{algorithmic}
\end{algorithm}
\section{Relation to other methods}\label{sec:others}
In this section we describe the relationship between the SS, the SURE and the GCV methods. In particular, we show that similarly to the SS method, the SURE and the GCV methods minimize an approximation of \(\text{PMSE}(\lambda)\) in \cref{eq:PMSEdefn}.
As discussed in \cref{sec:minPMSE}, the SS method approximates \cref{eq:PMSEdefn} by using the Picard parameter to split \(C(\lambda)\) in \cref{eq:TN} and by taking the expected value of only part of the sum, as written in \cref{eq:TnApprox}. In contrast, the SURE method approximates \(\text{PMSE}(\lambda)\) by taking the expected value of \(||n||^2\) and of the whole \(C(\lambda)\) without splitting it. The SURE method thus minimizes the function
\be\label{eq:SUREfun} \text{SURE}(\lambda) = \rho(\lambda) + ms^2 - 2s^2T(\lambda),\ee
where
\be\label{eq:Tredf}
T(\lambda) = \frac{\mathcal{E}(C(\lambda))}{s^2} = m-\text{rank}(A)+\sum_{k=r+1}^{r+q}\frac{\lambda^2}{\gamma_k^2+\lambda^2}.
\ee
To conclude, the SS method uses a more accurate approximation of \(C(\lambda)\) since the sums containing \(|\beta_k|^2\) in \cref{eq:TnApprox} capture at least part of the true, oscillatory behavior of \(\Re(\beta_k\overline{\nu}_k)\) in \cref{eq:TN}, in contrast to the SURE method, which replaces these terms with a constant.

Another popular method for determining \(\lambda\) is the GCV \cite{splines,GCV2}, which relies on the minimization of the function
\be\label{eq:GCVFun}
G(\lambda) = \frac{\rho(\lambda)}{(T(\lambda))^2},
\ee
where \(\rho(\lambda)\) is defined in \cref{eq:ResidualNorm} and \(T(\lambda)\) is defined in \cref{eq:Tredf}.
Despite the difference in forms between the SURE and the GCV functions, it can be shown that their minima are close to each other. Specifically, it is easy to show that
\be\label{eq:TwTr}
T(\lambda) = \text{trace}(I-H_{\lambda}),
\ee
where \(H_{\lambda} = A(A^*A+\lambda^2L^*L)^{-1}A^*\). To do so, we note that \cref{eq:TikhSolnUseless} implies \(H_{\lambda}b=Ax(\lambda)\) and therefore, recalling that \(\beta_k=u_k^*b\), we can rewrite \cref{eq:Ax_l} as
\be\label{eq:TwTrDeriv}H_{\lambda}b = \left(\sum_{k=1}^ru_ku_k^* + \sum_{k=r+1}^{r+q}\frac{\gamma_k^2}{\gamma_k^2+\lambda^2}u_ku_k^*\right)b.\ee
Since \cref{eq:TwTrDeriv} holds for any \(b\in\mathbb{R}^m\), we conclude that
\be\label{eq:HatMatInGSVD} H_{\lambda} = \sum_{k=1}^ru_ku_k^* + \sum_{k=r+1}^{r+q}\frac{\gamma_k^2}{\gamma_k^2+\lambda^2}u_ku_k^*.\ee
Since \(U\) is unitary, we can write
\[ I = UU^* = \sum_{k=1}^mu_ku_k^*,\]
and therefore,
\be\label{eq:HatMatDiff} I-H_{\lambda} = \sum_{k=r+1}^{r+q}\frac{\lambda^2}{\gamma_k^2+\lambda^2}u_ku_k^* + \sum_{k=r+q+1}^mu_ku_k^*.\ee
Expressing \cref{eq:HatMatDiff} in basis \(\{u_k\}\), we find that
\be\label{Eq:HatMatDiffInU}\begin{aligned} U^*\left(I-H_{\lambda}\right)U &= \sum_{k=r+1}^{r+q}\frac{\lambda^2}{\gamma_k^2+\lambda^2}\left(U^*u_k\right)\left(u_k^*U\right) + \sum_{k=r+q+1}^m\left(U^*u_k\right)\left(u_k^*U\right)\\ &= \sum_{k=r+1}^{r+q}\frac{\lambda^2}{\gamma_k^2+\lambda^2}e_ke_k^* + \sum_{k=r+q+1}^me_ke_k^*,\end{aligned}\ee
where \(\{e_k\}_{k=1}^m\) is the standard basis. Noting that \(e_ke_k^*=\text{diag}\{\underbrace{0,\ldots,0}_{k-1},1,\underbrace{0,\ldots,0}_{m-k}\}\), we obtain the desired result
\be\label{Eq:TrTwFinalDeriv} \text{trace}(I-H_{\lambda}) = \text{trace}\left(U^*(I-H_{\lambda})U\right) = \sum_{k=r+1}^{r+q}\frac{\lambda^2}{\gamma_k^2+\lambda^2} + m - \underbrace{(r + q)}_{=\text{rank}(A)} = T(\lambda).\ee
The function \(T(\lambda)\) is equivalent to the residual effective degrees of freedom used in regression analysis, see \cite[p. 63]{splines}. Therefore, the approximation
\be\label{eq:xTApprox}
\frac{\rho(\lambda^*)}{T(\lambda^*)} \approx s^2
\ee
holds when \(\lambda^*\) is the argument of the minimum of \cref{eq:GCVFun}, see \cite[sect. 4.7]{splines}, \cite[sect. 6.3]{LCurveAnal}.
Differentiating \cref{eq:SUREfun} and \cref{eq:GCVFun} twice, and using the approximation \cref{eq:xTApprox} and the fact that $T(\lambda),T'(\lambda), \rho'(\lambda)>0$ for all $\lambda\neq0$, it is easy to show that the local minima of both \(\text{SURE}(\lambda)\) in \cref{eq:SUREfun} and \(G(\lambda)\) in \cref{eq:GCVFun} satisfy
\be
\rho'(\lambda^*) = 2s^2T'(\lambda^*), \ \text{and}\ \rho''(\lambda^*)-2s^2T''(\lambda^*)> 0.
\ee
This result implies that both the GCV and the SURE methods have a local minimum at \(\lambda^*\). In practice, the global minimum coincides with this local minimum or is located very close to it. Therefore, we conclude that the GCV method is approximately equivalent to the SURE method and, consequently, that it too relies on replacing \(C(\lambda)\) in \cref{eq:TN} with its expected value \(s^2 T(\lambda)\) in \cref{eq:Tredf}. Thus, both the SURE and the GCV rely on an inferior approximation of \(C(\lambda)\),  compared with our SS method as explained above. For additional analysis of the relation between the SURE method and the GCV method see \cite{Li1985}.

\section{Numerical examples}\label{sec:numEx}

In this section we present the results of the Tikhonov regularization with the regularization parameter estimated using the following methods:
\begin{itemize}
    \item Regularization with the optimal regularization parameter minimizing the MSE \cref{eq:MSEdefn} (Tikh$_{OPT}$);
    \item The SS method minimizing the PMSE with \(k_0\) and \(s^2\), estimated by either \cref{alg:PicIndAlg} (SS$_{P3}$), \cref{alg:PicIndLilAlg} (SS$_{P2}$) or the algorithm that uses the Lilliefors test from \cite{PicInd1} (SS$_{PL}$);
    \item The SS method minimizing the MSE with \(k_0\) and \(s^2\), estimated by either \cref{alg:PicIndAlg} (SS$_{M3}$), \cref{alg:PicIndLilAlg} (SS$_{M2}$) or the algorithm that uses the Lilliefors test from \cite{PicInd1} (SS$_{ML}$);
    \item The DF method with \(k_0\) and \(s^2\), estimated using \cref{alg:PicIndAlg};
    \item The GCV method;
    \item The SURE method with \(s^2\), estimated using \cref{alg:PicIndAlg}.
\end{itemize}
To implement \cref{alg:PicIndAlg}, we set \(\varepsilon=5\times10^{-2}\) and \(h=\lceil m/50 \rceil\). The method however is generally robust to changes of up to an order of magnitude in these quantities.
We also note that occasionally, the algorithm for estimation of \(k_0\) given in \cite{PicInd1} returns indices outside the interval \([r+1,r+q]\). In these cases we set \(k_0\) to either \(r+1\) or \(r+q\), depending on which is closer to the estimated \(k_0\).

We test the above detailed methods on four test problems. The first problem is two-dimensional while the other three are one-dimensional problems taken from the \texttt{Regularization Tools} package \cite{RegTools}:
\begin{itemize}
    \item The \texttt{mri} test problem, an image-deblurring problem with the test image of size \(128\times128\) taken from \texttt{Matlab}'s image processing toolbox. The coefficient matrix \(A\) is chosen to be the separable Gaussian blur used in \cite[Sect. 5]{projRegGenForm}, which has full rank with condition number \(\text{cond}(A)=3.13\times10^{16}\). We set \(L=I\) to allow for a fair comparison to \cite{PicInd1} and note that this is a relatively large-scale problem.
    \item The \texttt{gravity} test problem for which we set \(L=D^{(1)}\), where \(D^{(m)}\) is the finite difference approximation of the \(m\)th derivative operator. This problem, while relatively small-scale, is severely rank-deficient with \(\text{rank}(A)=45\) for \(A\in\mathbb{R}^{1000\times1000}\).
    \item The \texttt{phillips} test problem, first introduced in \cite{Phillips1962}, for which we set \(L=D^{(2)}\). For this example, the coefficient matrix has full rank with \(\text{cond}(A)=2.64\times10^{10}\) for \(A\in\mathbb{R}^{1000\times1000}\).
    \item The \texttt{heat} test problem, for which we also set \(L=I\). This problem is only mildly rank-deficient with \(\text{rank}(A)=588\) for \(A\in\mathbb{R}^{1000\times1000}\).
\end{itemize}
For each problem, we add white Gaussian noise of zero mean and a variance of \(s^2=\alpha\max\{|b_{true}|^2\}\), where \(\alpha\in\{10^{-2}, 10^{-4}, 10^{-6}\}\). We thus present a total of twelve tests and for each test we generate 100 independent noise realizations.

To find the global minima of the functions associated with the above methods, we use the following algorithm. First, we evaluate a function on a sparse grid consisting of 1000 logarithmically spaced points between \(10^{-15}\) to \(10^{5}\). We then apply \texttt{Matlab}'s \texttt{fminbnd} solver with \texttt{TolX} of 1e-9 on \(11\) point intervals centered about each local minima found on the grid. Finally, we choose the value of \(\lambda\) corresponding to the global minimum found by the solver.

To assess the performance of each method, we use the mean-square deviation (MSD) defined as
\be\label{eq:MSD} \text{MSD}(\lambda) = \frac{||x_{true}-x(\lambda)||^2}{||x_{true}||^2}.\ee
The optimal solution is then defined as the one minimizing the MSD \cref{eq:MSD}, which also minimizes the MSE, since \(x_{true}\) does not depend on \(\lambda\).

\subsection{Results}\label{sec:results}
\begin{figure}[!p]
\centering
\includegraphics[width=\textwidth]{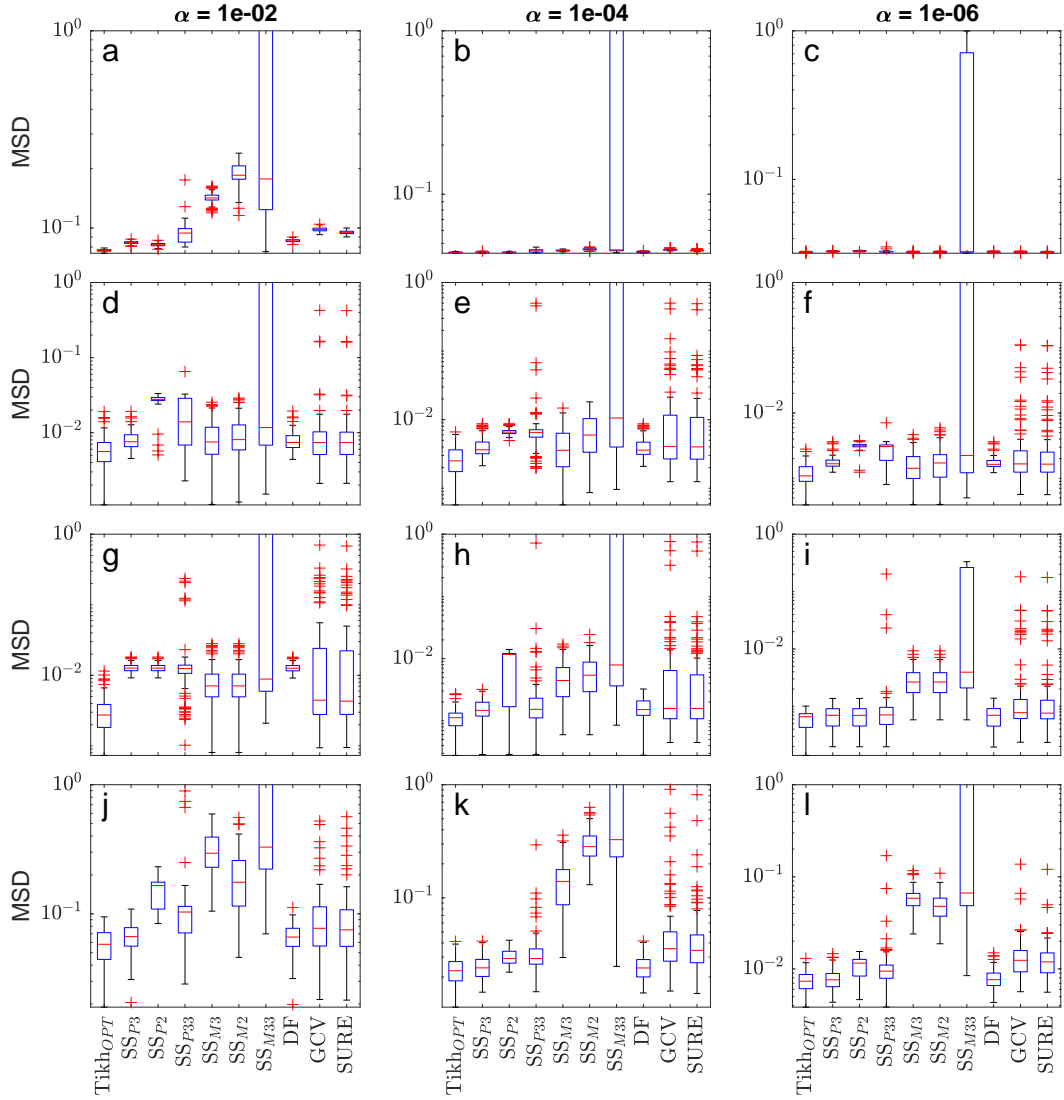} 
\caption{Boxplots of the MSD values of the reconstructions of the following examples, \emph{first row}: \texttt{mri}; \emph{second row}: \texttt{gravity}; \emph{third row}: \texttt{phillips}; \emph{fourth row}: \texttt{heat}. The noise levels presented are \emph{first column}: \(\alpha=10^{-2}\); \emph{second column}: \(\alpha=10^{-4}\); \emph{third column}: \(\alpha=10^{-6}\).}\label{fig:Boxplots}
\end{figure}

\begin{figure}[!p]
\centering
\includegraphics[width=\textwidth]{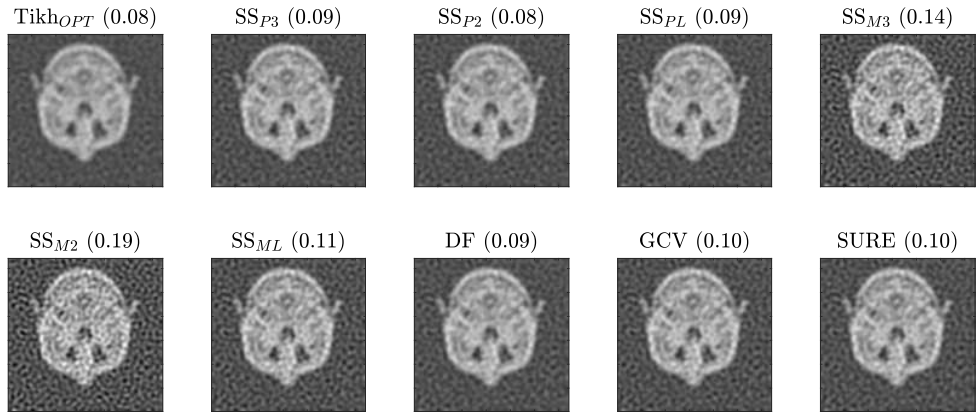} 
\caption{Reconstructed images for the \texttt{mri} test problem for noise level $\alpha=10^{-2}$. The MSD values \cref{eq:MSD} are listed in parentheses.}\label{fig:mri_rec_imgs}
\end{figure}

\begin{figure}[!p]
\centering
\includegraphics[width=\textwidth]{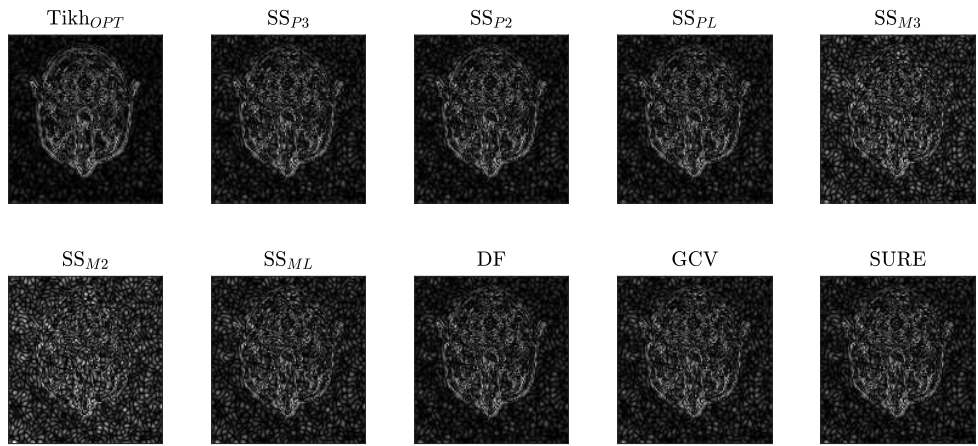} 
\caption{Error images $|x_{true}-x(\lambda)|$ (black = 0) for the \texttt{mri} test problem for noise level $\alpha=10^{-2}$ .}\label{fig:mri_err_imgs}
\end{figure}

In \cref{fig:Boxplots}, we present the results of our simulations by means of boxplots of the MSD values in log scale. Boxplots graphically depict the results by splitting them into quartiles so that each box spans the range between the first and third quartiles, termed the interquartile range (i.e., the middle 50\% of the data). The horizontal line in each box denotes the median and the error bars span 150\% of the interquartile range above the third quartile and below the first quartile. Any point outside this interval is denoted by '+' and  considered an outlier. We truncate the \(y\)-axis at \(\text{MSD}=1\) and present the number of truncated points for each method in \cref{tbl:FailRate}. For the \texttt{mri} test problem and noise level $\alpha=10^{-2}$, we also show the reconstructed images and the error images $|x_{true}-x(\lambda)|$ (with the absolute value applied pixel-wise) for all methods in \cref{fig:mri_rec_imgs} and \cref{fig:mri_err_imgs}, respectively.

From the results in \cref{fig:Boxplots}, we can make the following observations:
\begin{enumerate}
    \item The SS$_{P3}$ method performs consistently better than or very similarly to SS$_{P2}$, and both outperform SS$_{PL}$, which tends to fail in a significant percentage of cases, as can also be seen in \cref{tbl:FailRate}.
    \item For the \texttt{heat} problem with $L=I$ (\cref{fig:Boxplots}(j)-(l)), the method $SS_{ML}$ coincides with Algorithm 1 from \cite{PicInd1}, enabling a fair comparison to our methods. The MSE minimizing algorithms $SS_{M2}$ and $SS_{M3}$ are more accurate than $SS_{ML}$, however they are far from the optimal Tikhonov solution. Only by using both the minimization of the PMSE and our novel Picard estimation algorithm we obtain an almost optimal solution.
    \item The DF method performs almost identically to SS$_{P3}$ when both methods employ the same algorithm to estimate the Picard parameter.
    \item The SS methods that minimize the PMSE with the Picard parameter, found using any of the algorithms, performed similarly to or better than their counterparts minimizing the MSE in the majority of cases. The exceptions are the \texttt{gravity} example in \cref{fig:Boxplots}(d)-(f) and the \texttt{phillips} example with \(\alpha=10^{-2}\) in \cref{fig:Boxplots}(g). In the \texttt{gravity} example, the median of these methods is approximately the same, but the MSD values of the SS methods minimizing the MSE deviate more from the median compared to those of the SS methods minimizing the PMSE, both upwards and downwards. For the \texttt{phillips} example with \(\alpha=10^{-2}\), minimizing the MSE has a small advantage which is lost for lower noise levels. However, minimization of the PMSE also produces good results.
    \item In contrast to the consistent advantage of SS$_{P3}$ over the SS$_{P2}$ method, SS$_{M3}$ does not consistently outperform SS$_{M2}$. For example, in the \texttt{heat} test problem with noise levels \(\alpha=10^{-2}\) and \(\alpha=10^{-6}\) the SS$_{M2}$ is more accurate than the SS$_{M3}$ method. This is due to the fact that, contrary to the SS method minimizing the PMSE where the Picard parameter is only used to split the sum \(C(\lambda)\), in the SS method that minimizes the MSE the sums \(\widetilde \rho(\lambda)\) and \(\widetilde D(\lambda)\) are truncated. The truncation, even at the exact Picard parameter, decreases the accuracy of the approximation due to the lost terms. Hence, an overestimation of the Picard parameter will definitely decrease the accuracy of the PMSE approximation but may improve it for the MSE approximation by including a small number of additional terms in the minimized function, as long as division by small $\sigma_j\sigma_k$ does not cause instability.
    \item While the algorithm of \cite{PicInd1} for estimation of the Picard parameter is inferior to \cref{alg:PicIndLilAlg} and \cref{alg:PicIndAlg}, it performs better with the SS method minimizing the PMSE (SS$_{PL}$) than the one minimizing the MSE (SS$_{ML}$), and as shown in \cref{tbl:FailRate} is also more consistent. This difference is explained as above by the fact that truncation of the MSE approximation at the estimated Picard parameter leads to exclusion of a part of the high frequency information from the minimization target function.
    \item Both SS$_{P3}$ and the DF method outperformed the GCV and the SURE methods, which do not split the sums of the PMSE expansion as in \cref{eq:TnApprox}. In addition, the GCV and the SURE methods performed almost identically, as expected from the discussion in \cref{sec:others}.
\end{enumerate}

To conclude, our numerical examples clearly demonstrate the advantage of our new algorithm for estimation of the Picard parameter and of our approximation of the PMSE for both \(L=I\) and \(L\neq I\).

\section{Conclusions}

We generalized the approach taken in \cite{PicInd1} to estimate the regularization parameter for the general-form Tikhonov regularization. While the authors of \cite{PicInd1} approximately minimize the MSE, we show that such an approximation is numerically unstable in this generalized setting and propose to approximately minimize the PMSE instead. We develop two algorithms to stably approximate the PMSE, which we term the Series Splitting (SS) and Data Filtering (DF) methods, using the concept of the Picard parameter. While the two methods perform very similarly in the present framework, in which the SVD of the coefficient matrix $A$ is available, DF can be naturally generalized to large-scale problems, in which computing this SVD is prohibitive. This will be the subject of future work. We also present an algorithm similar to SS for a stable approximate minimization of the MSE, although the resulting approximation is less accurate than that for the PMSE due to the need to drop the numerically unstable terms. This algorithm can be used in the rare cases when the minimizers of the PMSE and MSE are far apart. The accuracy of all of these methods depends on the estimation of the Picard parameter, for which we proposed a novel algorithm based on an estimate of the variance of the noise.

Our methods were tested on multiple numerical examples and compared to the methods of \cite{PicInd1}, the GCV and the SURE. The numerical results indicate that, in contrast to other methods, the SS and the DF methods consistently produce near-optimal results for all test problems and noise realizations.

\section{Acknowledgements}
The authors would like to thank the referees and the editor for their constructive criticism and helpful suggestions to make this paper better.

\bibliographystyle{siamplain}
\bibliography{library}

\newgeometry{left=2.5cm,right=2.5cm}

\begin{table}[!p]
\centering
\caption{Number of points with MSD values larger than 1 out of a total of 100 points for each method. These points were truncated from the boxplots in \cref{fig:Boxplots}. The row abbreviations are consistent with those in \cref{fig:Boxplots}.}\label{tbl:FailRate}
\begin{tabular}{|c|ccccccccc|}
  \hline
 \diaghead{FigureMethod}{Fig.}{Met.}  & SS$_{P3}$ & SS$_{P2}$ & SS$_{PL}$ & SS$_{M3}$ & SS$_{M2}$ & SS$_{ML}$ & DF & GCV & SURE \\ \hline
  (a) & 0 & 0 & 0 & 0 & 0 & 42 & 0 & 0 & 0 \\
  (b) & 0 & 0 & 0 & 0 & 0 & 28 & 0 & 0 & 0 \\
  (c) & 0 & 0 & 0 & 0 & 0 & 24 & 0 & 0 & 0 \\
  (d) & 0 & 0 & 5 & 0 & 0 & 31 & 0 & 6 & 6 \\
  (e) & 0 & 0 & 1 & 0 & 0 & 32 & 0 & 3 & 3 \\
  (f) & 0 & 0 & 5 & 0 & 0 & 29 & 0 & 5 & 5 \\
  (g) & 0 & 0 & 2 & 0 & 1 & 25 & 0 & 8 & 6 \\
  (h) & 0 & 0 & 2 & 0 & 0 & 31 & 0 & 4 & 4 \\
  (i) & 0 & 0 & 1 & 0 & 0 & 24 & 0 & 1 & 0 \\
  (j) & 0 & 0 & 1 & 0 & 0 & 30 & 0 & 1 & 0 \\
  (k) & 0 & 0 & 3 & 0 & 0 & 25 & 0 & 0 & 0 \\
  (l) & 0 & 0 & 2 & 0 & 0 & 25 & 0 & 0 & 0 \\
  \hline
\end{tabular}
\end{table}
\restoregeometry

\end{document}